\documentclass[12pt,leqno]{amsart}
\usepackage[latin1]{inputenc}
\usepackage{color}
\usepackage{amssymb}
\usepackage{enumerate,tikz}
\usepackage{comment}
\usepackage[all]{xy}
\usepackage{hyperref,cleveref,graphics,mathrsfs}

\DeclareFontFamily{U}{matha}{\hyphenchar\font45}
\DeclareFontShape{U}{matha}{m}{n}{
      <5> <6> <7> <8> <9> <10> gen * matha
      <10.95> matha10 <12> <14.4> <17.28> <20.74> <24.88> matha12
      }{}
\DeclareSymbolFont{matha}{U}{matha}{m}{n}
\DeclareFontSubstitution{U}{matha}{m}{n}

\DeclareFontFamily{U}{mathx}{\hyphenchar\font45}
\DeclareFontShape{U}{mathx}{m}{n}{
      <5> <6> <7> <8> <9> <10>
      <10.95> <12> <14.4> <17.28> <20.74> <24.88>
      mathx10
      }{}
\DeclareSymbolFont{mathx}{U}{mathx}{m}{n}
\DeclareFontSubstitution{U}{mathx}{m}{n}

\DeclareMathDelimiter{\vvvert}{0}{matha}{"7E}{mathx}{"17}

\usetikzlibrary{matrix}

\usepackage{lmodern} 
\usepackage{amsmath} 
\usepackage{mathtools} 
\usepackage{commath} 
\usepackage{amsfonts} 
\usepackage{dsfont} 

\usepackage{graphicx} 
\usepackage{amsthm} 
\numberwithin{equation}{section} 

\newtheorem{thm}{Theorem}[section]

\newtheorem{lem}[thm]{Lemma}

\theoremstyle{definition}
\newtheorem{defn}[thm]{Definition}

\newtheorem{rem}[thm]{Remark}

\newtheorem{example}[thm]{Example}
\numberwithin{equation}{section}

\newcommand{\eps}{\varepsilon}

\def\epsilon{\varepsilon}

\usepackage{bbm}
\newcommand{\vertiii}[1]{{\left\vert\kern-0.25ex\left\vert\kern-0.25ex\left\vert #1 
    \right\vert\kern-0.25ex\right\vert\kern-0.25ex\right\vert}}

\newcommand{\triple}[1]{{\left\vert\kern-0.25ex\left\vert\kern-0.25ex\left\vert #1 
    \right\vert\kern-0.25ex\right\vert\kern-0.25ex\right\vert}}

\newcommand{\N}{\mathbb{N}}

\newcommand{\R}{\mathbb{R}}
\newcommand{\C}{\mathbb{C}}

\newcommand{\Lbb}{\mathbb{L}}

\newcommand{\Ccal}{\mathcal{C}}
\newcommand{\Fcal}{\mathcal{F}}

\newcommand{\Hcal}{\mathcal{H}}

\newcommand{\Gcal}{\mathcal{G}}

\newcommand{\Dcal}{\mathcal{D}}

\def\Lip{\operatorname{Lip}}
\def\Span{\mathrm{span}}

\newcommand{\spn}{\mathrm{span}}
\newcommand{\fbl}{\mathrm{FBL}}
\newcommand{\fbp}{{\mathrm{FBL}}^{(p)}}

\DeclareMathOperator{\FVL}{FVL}
\DeclareMathOperator{\FBLi}{FBL^{(\infty)}}

\DeclareMathOperator*{\latt}{\text{lat}}

\newcommand{\nFBL}[2]{\norm{#1}_{\mathrm{FBL}[#2]}}
\newcommand{\nFBLE}[1]{\norm{#1}_{\mathrm{FBL}[E]}}

\usepackage{bbm}

\title{Free objects in Banach space theory}

\subjclass[2020]{46B42, 46B10, 46G20, 46E15, 46B20} 

\author[E. Garc\'ia-S\'anchez]{Enrique Garc\'ia-S\'anchez}
\address{Instituto de Ciencias Matem\'aticas (CSIC-UAM-UC3M-UCM)\\
Consejo Superior de Investigaciones Cient\'ificas\\
C/ Nicol\'as Cabrera, 13--15, Campus de Cantoblanco UAM\\
28049 Madrid, Spain.}
\email{enrique.garcia@icmat.es}

\author[D. de Hevia]{David de Hevia}
\address{Instituto de Ciencias Matem\'aticas (CSIC-UAM-UC3M-UCM)\\
Consejo Superior de Investigaciones Cient\'ificas\\
C/ Nicol\'as Cabrera, 13--15, Campus de Cantoblanco UAM\\
28049 Madrid, Spain.}
\email{david.dehevia@icmat.es}

\author[P.~Tradacete]{Pedro Tradacete}
\address{Instituto de Ciencias Matem\'aticas (CSIC-UAM-UC3M-UCM)\\
Consejo Superior de Investigaciones Cient\'ificas\\
C/ Nicol\'as Cabrera, 13--15, Campus de Cantoblanco UAM\\
28049 Madrid, Spain.}
\email{pedro.tradacete@icmat.es}

\date{\today}

\keywords{Free object; Lipschitz free space; free Banach lattice}

\thanks{Research partially supported by grants PID2020-116398GB-I00 and CEX2019-000904-S funded by MCIN/AEI/10.13039/501100011033. E. Garc\'ia-S\'anchez was partially supported by FPI grant CEX2019-000904-S-21-3 funded by MCIN/AEI/10.13039/501100011033 and by ``European Union NextGenerationEU/PRT''. D. de Hevia benefited from an FPU Grant FPU20/03334 from the Ministerio de Universidades. P. Tradacete is also supported by a 2022 Leonardo Grant for Researchers and Cultural Creators, BBVA Foundation.}

\begin{document}

\begin{abstract}
    We survey recent progress on three relevant instances of free objects related to Banach spaces: Lipschitz free spaces generated by metric spaces, holomorphic free spaces generated by open sets and free Banach lattices generated by Banach spaces. Our emphasis will be on the parallelisms among these developments.
\end{abstract}

\maketitle

\section{Introduction}

The aim of this note is to review recent developments concerning free objects arising in different aspects of analysis and which are deeply connected to Banach space theory. The first of these notions that we will revisit here are \emph{Lipschitz free Banach spaces}, which associate a Banach space to a metric space. This construction has been the object of intense research in the last 20 years and has been useful in the study of the non-linear geometry of Banach spaces, with ramifications and applications in metric embeddings, geometric group theory... (cf. \cite{BL,CKN,Naor}). Secondly, we will review what we will call \emph{holomorphic free Banach spaces}, which originate in a construction due to J. Mujica \cite{Mujica91}, allowing to linearize holomorphic maps from an open set to a Banach space, and which, although not so developed, somehow parallel the theory of Lipschitz free spaces. Finally, we will focus on \emph{free Banach lattices}, which are a recent addition due to A. Avil\'es, J. Rodr\'iguez and P. Tradacete \cite{ART} and provide a new tool for studying the interactions between lattice structures and Banach spaces.  \\

All these examples can be framed within the notion of the \emph{free object} in a certain category generated by an object from a (typically larger) category, which has been central in several developments in different areas of mathematics. Probably, the most significant and better known examples are related to universal algebra, and include in particular free groups, free modules, free algebras or free lattices. \\

The language of Category theory is particularly convenient to deal with this notion. The idea of a free object can be somehow seen as a generalization of the notion of basis of a vector space, based on the property that every linear map is determined by the values on the basis, and has its origins in the work of S. MacLane \cite{Maclane}, J. R. Isbell \cite{Isbell} and Z. Semadeni \cite{Semadeni}, among others. The extension to the categories setting can be expressed as follows: Suppose we have categories $\mathcal C_1$ and $\mathcal C_2$, with a faithful functor $f:\mathcal C_2\rightarrow \mathcal C_1$; typically, we will consider $\mathcal C_1$ as a category which is ``larger'' than $\mathcal C_2$ (for instance, we could take $\mathcal C_1$ as the category of metric spaces with Lipschitz functions and $\mathcal C_2$ that of Banach spaces with bounded linear operators) and consider $f$ as the ``forgetful'' functor (which in the same example would take a Banach space $(X,\|\cdot\|)$ to the metric space $(X,d(\cdot,\cdot))$, with $d(x,y)=\|x-y\|$, and a bounded linear operator to the corresponding Lipschitz map); now, given an object $O$ in $\mathcal C_1$, the \textit{free object in $\mathcal C_2$ generated by $O$}, is an object $F(O)$ in $\mathcal C_2$ together with a morphism $i:O\rightarrow f(F(O))$ satisfying the following universal property: for any object $X$ in $\mathcal C_2$ and any morphism $\varphi:O\rightarrow f(X)$ in $\mathcal C_1$, there is a unique morphism $F(\varphi):F(O)\rightarrow X$ in $\mathcal C_2$ such that $f(F(\varphi))\circ i=\varphi$. It will be customary to plot this in the shape of the following commuting diagram: 
$$
\xymatrix{f(F(O))\ar@{-->}[rd]^{f(F(\varphi))}&\\
 O\ar[r]^{\varphi}\ar^{i}[u]& f(X)}
$$

Free objects have also played a fundamental role in topology, although these can be found in the literature under somewhat different terminologies. For instance, recall that the \textit{Stone-\v{C}ech compactification} of a topological space $S$ is a compact Hausdorff space $\beta S$ together with a continuous map $i_S:S\rightarrow \beta S$, so that for every compact Hausdorff space $K$ and every continuous map $f:S\rightarrow K$ there is a unique continuous map $\beta f:\beta S\rightarrow K$ making the following diagram commute:
$$
\xymatrix{\beta S\ar@{-->}[rd]^{\beta f}&\\
 S\ar[r]^{f}\ar^{i_S}[u]& K}
$$
Thus, we could consider $\beta S$ as the \textit{free compact Hausdorff space generated by} $S$.\\

A similar scheme can be used to describe different notions of completions or closures of a given space. These include in particular the completion of a metric, or uniform, space (i.e. adding limits of Cauchy sequences) or the order completion of lattices or boolean algebras (by adding infima and suprema of order bounded sets). \\

The study of free objects in Banach space theory can be considered as a more recent development. However, we will see next how the second dual of a Banach space can be framed in this picture.
Consider the category of \textit{dual} Banach spaces together with \textit{adjoint operators} (or equivalently, weak$^*$ continuous linear maps). Given a Banach space $E$, let $J_E:E\rightarrow E^{**}$ denote the canonical embedding; it is clear that every bounded linear operator $T:E\rightarrow X^*$ can be uniquely extended to a dual operator $\hat T:E^{**}\rightarrow X^*$, given by $\hat{T}=(T^*\circ J_X)^*$, in such a way that the following diagram commutes:
$$
\xymatrix{E^{**}\ar@{-->}[rd]^{\hat{T}}&\\
     E\ar[r]^{T}\ar^{J_E}[u]& X^*}
$$
Thus, we can consider $E^{**}$ (together with the embedding $J_E$) as the \textit{free dual Banach space generated by} $E$.\\

Note this construction however, is actually closer to that of a completion: if $E$ is a Banach space and we consider its unit ball $B_{E}$ equipped with the weak topology (the one induced by the functionals in $E^*$), then $B_{E^{**}}$ can be considered as the completion with respect to that topology (which actually makes the ball into a compact set by Banach-Alaoglu's theorem).\\


The paper is organized as follows. Section \ref{sect: Lipschitz free} is devoted to Lipschitz free Banach spaces generated by a pointed metric space. We will recall first the different ways to construct this space and then focus on the interaction between properties of the Lipschitz free space $\mathcal F(M)$ and the underlying metric space $M$. In Section \ref{sec:holomorfo} we introduce the notion of holomorphic free Banach space $\mathcal G^\infty(U)$ generated by an open subset $U$ of complex Banach space. After showing the construction of this object, we present an overview about its main properties. Finally, Section \ref{sec:fbl} focuses on the free Banach lattice generated by a Banach space, following a similar structure as in the previous sections.\\

Finally, note that an abstract approach to free objects can be useful for identifying general principles that have different manifestations in each of the aforementioned constructions. This approach has been pursued for instance in \cite{AMR3} exhibiting the relation between local complementation and extension properties in $\mathcal F(M)$ and $\fbl[E]$. A general construction for linearizing functions with values in locally convex spaces has also been considered in \cite{CZ} which covers in particular the case of $\mathcal G^\infty(U)$.\\

We refer the reader to the monographs \cite{BanachSpaceTheory, LT1} for Banach space theory, \cite{LT2, M-N} for Banach lattices and \cite{BL, Weaver} for the theory of Lipschitz functions.

\section{Lipschitz free Banach space}\label{sect: Lipschitz free}

In the recent years, especially after the appearance of \cite{Weaver}, \textit{Lipschitz free spaces} (also known as \textit{Arens-Eells spaces} or \textit{transportation cost spaces}) have attracted considerable attention from researchers interested in metric geometry and the non-linear theory of Banach spaces (see, for instance, the survey paper \cite{Godefroy}, and also \cite[1.6]{OO19} for historical remarks about this construction). These spaces can be defined by the following universal property: Given a metric space $M$ with a distinguished point $0$, $\mathcal F(M)$ is a Banach space together with an isometric map $\delta:M\rightarrow \mathcal F(M)$ with the property that for every Banach space $X$ and every Lipschitz map $f:M\rightarrow X$ with $f(0)=0$, there is a unique linear operator $\widehat f:\mathcal F(M)\rightarrow X$ making  the following diagram commute:
$$
\xymatrix{\mathcal F(M)\ar@{-->}^{\widehat f}[rd]&\\
     M\ar[r]^{f}\ar^{\delta}[u]& X}
$$
Moreover, the operator norm of the linear extension coincides with the Lipschitz norm of the function $f$: $\|\widehat f\|=\|f\|_{Lip}.$ A very fruitful line of research is devoted to analyzing the interplay between Banach space properties of $\mathcal F(M)$ versus metric properties of $M$ (see for example \cite{AGPP, AP, cutDou, CDW, god, GK, GO}). 

\subsection{Construction}

Let $(M,d)$ be a metric space with a distinguished point denoted by $0$ and let $\Lip_0(M)$ denote the space of real-valued Lipschitz functions that map $0\in M$ to $0\in\mathbb R$. We can equip this space with the norm $$\|f\|_{\Lip}=\sup\Big\{\frac{|f(x)-f(y)|}{d(x,y)}:x\neq y\in M\Big\},$$ for $f\in \Lip_0(M)$, so that $\Lip_0(M)$ becomes a Banach space.\\

For $x\in M$, let $\delta_x\in \Lip_0(M)^*$ be the evaluation functional given by $\delta_x(f)=f(x)$ for $f\in \Lip_0(M)$. Let us denote by $\mathcal F(M)$ the closed linear span of $\{\delta_x:x\in M\}$ with the dual norm in $\Lip_0(M)^*$. Note that the map $x\mapsto \delta_x$ becomes an isometry, as we have $\|\delta_x-\delta_y\|=d(x,y)$ for any $x,y\in M$.\\

The space $\mathcal F(M)$ together with $\delta_M:M\rightarrow \mathcal F(M)$ is the \textit{Lipschitz free Banach space} and it is uniquely characterized by the following universal property:\\

For every Banach space $X$ and every Lipschitz map $f:M\rightarrow X$ with $f(0)=0$, there is a unique linear operator $\widehat f:\mathcal F(M)\rightarrow X$ making  the following diagram commute:
$$
\xymatrix{\mathcal F(M)\ar@{-->}^{\widehat f}[rd]&\\
     M\ar[r]^{f}\ar^{\delta}[u]& X}
$$

Let us briefly recall the proof of this fact. Fix a Banach space $X$ and a Lipschitz map $f:M\rightarrow X$ mapping $0$ to $0$. Extend linearly $f$ to $\Span\{\delta_x: x\in M\}$ and denote this extension by $\widehat{f}$. It is enough to check that $\|\widehat{f}\|=\|f\|_{\Lip}$. To this end, pick  $a\in\Span\{\delta_x: x\in M\}$. Then $\|\widehat{f}(a)\|_X=x^*(\widehat{f}(a))$ for some $x^*\in B_{X^*}$. However, $x^*\circ f$ belongs to $\Lip_0(M)$ and $\|x^*\circ f\|_{\Lip}\leq \|f\|_{\Lip}$, so $\|a\|\|f\|_{\Lip}\geq \|\widehat{f}(a)\|_X$ as claimed. Hence $\widehat{f}$ can be uniquely extended to $\mathcal F(M)$. As a consequence we can identify $\mathcal F(M)^* = \Lip_0(M)$. \\

The above universal property can also be extended in a natural way to another one with a more functorial character. Suppose $(M,d)$, $(N,d')$ are pointed metric spaces, then every Lipschitz map $f:M\rightarrow N$ mapping $0$ to $0$ admits a linear extension $\overline f$ making the following diagram commutative:

$$
\xymatrix{\mathcal F(M)\ar@{-->}[r]^{\overline{f}}& \mathcal F(N)\\
     M\ar[r]^{f}\ar^{\delta_M}[u]& N \ar^{\delta_N}[u]}
$$

It is not hard to see that the map $f$ is bi-Lipschitz if and only if $\overline f$ is an isomorphic embedding, that $f$ is a Lipschitz retract if and only if $\overline f$ is a linear projection, or that $f$ has dense range if and only if $\overline f$ has. The question about injectivity is less clear and has been pursued in \cite{GPP}.\\

In connection with this, an important consequence of McShane's extension theorem for Lipschitz functions \cite[Theorem 1.5.6]{Weaver} is the following: Suppose $N$ is a subspace of a metric space $(M,d)$, let $\iota:N\hookrightarrow M$ denote the formal inclusion, which is clearly Lipschitz, then $\overline{\iota}$ defines a linear isometric embedding of $\mathcal F(N)$ as a subspace of $\mathcal F(M)$. As a consequence, we get that if $N$ is bi-Lipschitz equivalent to a subset of $M$, then $\mathcal F(N)$ is linearly isomorphic to a subspace of $\mathcal F(M)$.\\

There is another intrinsic way to define Lipschitz free spaces by means of the notion of the Arens-Eells norm or transportation cost spaces, which are closely related to Wasserstein-1 spaces \cite{OO19}. For further details and motivation we refer to \cite[Section 2.3]{Weaver}, but we include here the main steps in this idea: For $a\in \Span\{\delta_x: x\in M\}$, let us consider
$$\vvvert a \vvvert=\inf\Big\{\sum_{j=1}^{n}|\alpha_j| d(y_j,z_j): a=\sum_{j=1}^{n}\alpha_j(\delta_{y_j}-\delta_{z_j})\Big\}.$$
On the one hand, it can be easily checked that $\vvvert\cdot\vvvert$ is the largest seminorm $\rho$ on $\Span\{\delta_x: x\in M\}$ satisfying $\rho(\delta_x-\delta_y)\leq d(x,y)$ for every $x,y\in M$. As the standard norm $\|\cdot\|_{\mathcal F(M)}$ on $\mathcal F(M)$ satisfies this condition, it follows that $\|\cdot\|_{\mathcal F(M)}\leq\vvvert\cdot\vvvert$. In particular, $\vvvert\cdot\vvvert$ is actually a norm and $\vvvert\delta_x-\delta_y\vvvert=d(x,y)$ for $x,y\in M$.\\

On the other hand, if we consider now the Banach space $$X=\overline{\Span\{\delta_x: x\in M\}}^{\vvvert\cdot\vvvert}$$ and take the map $f:M\rightarrow X$ given by $f(x)=\delta_x$, then $f$ is an isometric embedding, and by the universal property extends to a linear map $\widehat{f}:\mathcal F(M)\rightarrow X$ which is $1$-Lipschitz. Hence, it follows that $X=\mathcal F(M)$ with $\vvvert\cdot\vvvert=\|\cdot\|_{\mathcal F(M)}$, which is the claim.

\subsection{Banach space properties of $\mathcal F(M)$}\label{ss:propertieslipfree}

As a predual of the space of Lipchitz functions, $\mathcal F(M)$ has a certain $\ell_1$ behavior. It is well known that if $M\subset\mathbb R$ is a set of measure zero or if $M$ is a separable ultrametric space, then $\mathcal F(M)$ is isomorphic to $\ell_1$, see \cite{god} and \cite{cutDou}.\\ 

In \cite{CDW} the authors showed that for every infinite metric space $M$, $\ell_1$ always embeds as a complemented subspace of $\mathcal F(M)$. As a consequence, we get that $\mathcal F(M)$ is never isomorphic to a complemented subspace of a $C(K)$ space, or that $\mathcal F(M)^*$ is not weakly sequentially complete (see definition below). Further results about the relation with $\ell_1$ can be found in \cite{APP,OO20}.\\

Recall that a Banach space $X$ has the \textit{approximation property} (AP) if the identity $id_X$ of $X$ is in the closure of the finite rank operators on $X$ for the topology of uniform convergence on compact sets. The \textit{$\lambda$-bounded approximation property} ($\lambda$-BAP) means that there are finite rank operators with norm less than $\lambda$ that can be chosen in the approximation, and the (1-BAP) is called the \textit{metric approximation property} (MAP).\\

In \cite{GO} it was shown that for a separable metric space  $M$ and $\lambda\geq1$, the Lipschitz free space $\mathcal F(M)$ has the $\lambda$-BAP if and only if for any Banach space $Y$ and any $f\in \Lip(M,Y^{**})$, there is a net $(f_\alpha)\subset \Lip(M,Y)$ which converges to $f$ in the pointwise-weak$^*$ topology with $\|f_\alpha\|\leq\lambda\|f\|$. Note that in the particular case when $M$ is a Banach space, \cite[Theorem 5.3]{GK} already shows that $\mathcal F(X)$ has the $\lambda$-BAP if and only if $X$ has the $\lambda$-BAP. More recently, it was shown in \cite{LP} that for every doubling metric space $M$, the space $\mathcal F(M)$ has BAP, and actually $\mathcal F(M)$ has MAP for every compact convex $M\subset \mathbb R^n$ (see \cite{PS}). Going even one step further, it was shown in \cite{HP} that $\mathcal F(\mathbb R^n)$ and $\mathcal F(\ell_1)$ actually have a Schauder basis. On the opposite direction, in \cite{HLP} an example is given of a compact metric space homeomorphic to the Cantor set whose Lipschitz free space fails the AP.\\

As noted above,  $\mathcal F(M)^* = \Lip_0(M)$. It was asked in \cite{Weaver} whether the space $\Lip_0(M)$ has a unique predual. This has been settled affirmatively in \cite{Weaver18} for the cases when $M$ has finite diameter or when $M$ is a Banach space.\\

In connection with duality, it seems to be a difficult problem to identify the extreme points of the unit ball of $\mathcal F(M)$. Natural candidates for such points are the elements of the form $$m_{x,y}=\frac{\delta_x-\delta_y}{d(x,y)},$$ usually called \textit{(elementary) molecules}. It was proved in \cite{GPPR18} that the sets of preserved extreme points and denting points of the unit ball of $\mathcal F(M)$ coincide. Also, in \cite{GPR18} the authors show that a molecule $m_{x,y}$ is a strongly exposed point of the unit ball of $\mathcal F(M)$ if and only if
$$
d(x,z)+d(y,z)-d(x,y)\geq C\min\{d(x,z),d(y,z)\}
$$
for some $C>0$ and every $z\in M$. Here, recall that a norm one element $x$ of a Banach space $X$ is called a \textit{strongly exposed point of $B_X$} if there is $x^*\in X^*$ with $\|x^*\|=1$ such that every sequence $(x_n)\subset B_X$ with $x^*(x_n)\rightarrow x^*(x)$ must be norm convergent to $x$. More recently, it has been shown in \cite{Aliaga} that if $M$ is a proper metric space (where bounded closed subsets are compact), the extreme points are precisely the molecules $m_{x,y}$ satisfying $d(x,z)+d(y,z)>d(x,y)$ for every $z\in M\backslash\{x,y\}$. It follows that all extreme points are preserved extreme points and also exposed points. See also \cite{AG, APPP} for further related results on extreme points. Apparently, it is still unknown whether all extreme points of the unit ball of $\mathcal F(M)$ are molecules for an arbitrary metric space $M$.\\

A classical geometrical property of a Banach space is the \textit{Daugavet property}, which can be defined by the equation $\|T+I\|=1+\|T\|$ holding for every rank-one operator $T$. The Daugavet property of the space $\mathcal F(M)$ has been characterized in \cite{GPR18} precisely as the property of the metric space $M$ to be a length space; this complements earlier results from \cite{IKW}. Recall that a metric space $M$ is called a \textit{length space} if, for each pair of points $x,y\in M$, $d(x,y)$ coincides with the infimum of the length of rectifiable curves joining $x$ and $y$ (when this infimum is attained we say that $M$ is a \textit{geodesic space}).\\

A Banach space has the \textit{Schur property} when every weakly convergent sequence must be norm-convergent. It was shown in \cite{HLP} that the free space over a countable compact metric space has the Schur property. In connection with this, a more general sufficient condition for $\mathcal F(M)$ to have the Schur property has been provided in \cite{Pet} strongly based on the work of N. Kalton \cite{K:04}: Let $$\text{lip}_0(M)=\Big\{f\in\Lip_0(M):\lim_{\varepsilon\rightarrow0}\sup_{0<d(x,y)<\varepsilon}\frac{|f(x)-f(y)|}{d(x,y)}=0\Big\};$$ if $\text{lip}_0(M)$ is $1$-norming for $\mathcal F(M)$, then $\mathcal F(M)$ has the Schur property. More recently, in \cite{AGPP}, the Schur property of $\mathcal F(M)$ has been characterized by the fact that (the completion of) $M$ is purely 1-unrectifiable, and this is also equivalent to the Radon-Nikod\'ym property of $\mathcal F(M)$. \\

Recall that a Banach space $X$ is called \textit{weakly sequentially complete} if every weakly Cauchy sequence in $X$ is weakly convergent. Also, a Banach space $X$ is \textit{super-reflexive} if no non-reflexive Banach space is finitely representable in $X$. If $M$ is a compact subset of a superreflexive Banach space, then the space $\mathcal F(M)$ is weakly sequentially complete \cite{KP}. Also, $\mathcal F(X)$ is weakly sequentially complete whenever $X$ is superreflexive \cite{APP21}. \\

It is clear that if $M_1$ and $M_2$ are Lipschitz isomorphic metric spaces, then the corresponding spaces $\mathcal F(M_1)$ and $\mathcal F(M_2)$ are linearly isomorphic. The converse of this is false as shown in \cite{DF}, where the authors showed that if $K$ is any infinite metric compact space, then $\mathcal F(C(K))$ and $\mathcal F(c_0)$ are linearly isomorphic.\\

The case when the underlying metric space $M$ is already a Banach space is particularly relevant. If $X$ is a Banach space, by construction there is a unique linear map $\beta_X:\mathcal F(X) \rightarrow X$ (usually called \emph{barycenter map}) extending the identity on $X$, that is $\beta_X \circ \delta_X=id_X$. It is natural to wonder under which conditions $\beta_X$ has a linear right inverse. This is equivalent to asking whether the Banach space $X$ linearly embeds as a subspace of $\mathcal F(X)$ so that $\beta_X$ defines a projection onto this subspace. This property was considered in the seminal paper by G. Godefroy and N. Kalton \cite{GK}, where it was shown that every separable Banach space and every space of the form $\mathcal F(M)$ have this property.\\

The work \cite{GK} also provides an important application of Lipschitz free spaces to the non-linear theory of Banach spaces: Namely, if a separable Banach space $X$ isometrically embeds into another Banach space $Y$, then there exists an isometric linear embedding, while this does
not hold in the non-separable setting.\\

Another interesting feature of $\mathcal F(X)$ for $X$ a Banach space proved in \cite{Kauf} is that $\mathcal F(X)$ is linearly isomorphic to the infinite $\ell_1$-sum of $\mathcal F(X)$. As a consequence, it is observed also in \cite{Kauf} that $\mathcal F(X)$ is linearly isomorphic to $\mathcal F(B_X)$.\\

Finally, let us note that a quasi-normed version of Lipschitz free spaces has been considered in \cite{AK}. Recall, $(X,\rho)$ is a \textit{quasi-metric space} if $\rho:X\times X\rightarrow [0,\infty)$ is a symmetric function such that $\rho(x,y)=0$ if and only if $x=y$ and there exists $k\geq1$ such that $\rho(x,y)\leq k(\rho(x,z)+\rho(y,z))$ for $x,y,z\in X$. In this direction, a vector space $X$, together with a map $\|\cdot\|:X\rightarrow [0,\infty)$ such that $\rho(x,y)=\|x-y\|$ is a quasi-metric and $\|\lambda x\|=|\lambda|\|x\|$ is called a \textit{quasi-normed space}. The Aocki-Rolewicz theorem states that every quasi-Banach space admits an equivalent $p$-norm $\vvvert\cdot\vvvert$ (that satisfies $\vvvert x+y\vvvert^p\leq \vvvert x\vvvert^p+\vvvert y\vvvert^p$). Now, for a pointed $p$-metric space $(X,\rho)$, one can construct the Lipschitz free $p$-space over $(X,\rho)$ in a similar spirit as that of $\mathcal F(M)$ so that the latter becomes the case $p=1$. This has been further explored in \cite{AACD, AACD20, AACD21}. Along similar lines free spaces over assymetric spaces (those where the ``distance'' function is not necessarily symmetric) have been considered in \cite{DSV}.\\

\section{Holomorphic free Banach space}\label{sec:holomorfo}

In the same spirit as Lipschitz free Banach spaces, we can construct what could be called \emph{holomorphic free Banach spaces}. This object is a complex Banach space which allows us to extend bounded holomorphic functions defined on a fixed open subset of a (complex) Banach space to continuous $\mathbb{C}$-linear operators.
Given complex Banach spaces $E$ and $F$ and an open subset $U$ of $E$, we say that a function $f:U\to F$ is \textit{holomorphic} if it is Fr\'echet differentiable at every $x\in U$, that is, for every $x\in U$ there exists a continuous $\mathbb{C}$-linear mapping $Df(x):E\to F$ such that 
\begin{equation}\label{eq:frechet-derivative}
 \lim_{h\to 0}\frac{f(x+h)-f(x)-Df(x)(h)}{\|h\|}=0.   
\end{equation}

For simplicity, throughout this section Banach spaces will always be defined over the field of complex numbers and operators on them will be $\mathbb{C}$-linear.

\subsection{Construction.} Given an open subset $U$ of a Banach space $E$, there is a Banach space $\mathcal{G}^\infty(U)$ together with a bounded holomorphic function $\delta_U:U\to\mathcal{G}^\infty(U)$ with the following universal property: for every Banach space $F$ and every bounded holomorphic function $f:U\to F$, there exists a unique operator $\hat{f}:\mathcal{G}^\infty(U)\to F$ such that $\hat{f}\circ \delta_U=f$ and $\norm[0]{\hat{f}}=\norm{f}_{\infty}$, where $\norm{f}_{\infty}=\sup_{x\in U}\abs{f(x)}$. In other words, $\hat{f}$ makes commutative the following diagram:
$$
\xymatrix{\mathcal{G}^\infty(U)\ar@{-->}[rd]^{\hat{f}}\\
     U\ar[r]^{f}\ar^{\delta_U}[u]& F }
$$

The existence of such an object was proved by Mujica in \cite[Theorem 2.1]{Mujica91}, inspired by the Dixmier-Ng theorem (see also \cite[p. 417]{Dineen}). The uniqueness of $\mathcal{G}^\infty(U)$ (up to an isometric isomorphism) follows directly from its definition. In  \cite{Mujica91}, the author showed that
$$
\mathcal{G}^\infty(U)=\{u\in\mathcal{H}^\infty(U)^*\::\:\left.u\right|_{(B_{\mathcal{H}^\infty(U),\tau_c})} \text{ is continuous}\},
$$
where
$$\mathcal{H}^\infty(U)=\{f:U\rightarrow \C: f \text{ is holomorphic and bounded}\},$$
$\tau_c$ represents the compact-open topology, and
$$\fullfunction{\delta_U}{U}{\mathcal{G}^\infty(U)}{x}{\delta_x}$$
is given by $\delta_x(f)=f(x)$, $f\in\mathcal{H}^\infty(U)$.

\begin{rem}
Note that from the universal property (setting $F=\C$) it follows that the holomorphic free Banach space $\mathcal{G}^\infty (U)$ is an isometric predual of $\mathcal{H}^\infty(U)$. A reasonable question that arises is whether for every open subset $U$ of a Banach space $E$ the space $\mathcal{H}^\infty(U)$ has a unique isometric predual (compare to the results mentioned above in Section \ref{sect: Lipschitz free} about the predual of $\Lip_0(M)$). T. Ando \cite{Ando} and P. Wojtaszczyk \cite{W} independently proved that $\mathcal{H}^\infty(\mathbb{D})$ has a unique predual. This has been generalized for other planar domains in \cite{RY}. However, it is not known whether $\mathcal{H}^\infty(\mathbb{D}^n)$ has a unique predual for $n>2$.
\end{rem}

A remarkable property of these spaces is that they always satisfy a sort of \textit{lifting property}. Given an open bounded subset $U$ of a Banach space $E$, we can consider the restriction to $U$ of the identity mapping on $E$, denoted by $\iota:U\to E$, which is a bounded holomorphic function. By the universal property of $\mathcal{G}^\infty(U)$, there exists an operator $\beta$ which makes commutative the diagram below
$$
\xymatrix{\mathcal{G}^\infty(U)\ar@{-->}[rd]^{\beta}\\
     U\ar[r]^{\iota}\ar^{\delta_U}[u]& E }
$$
It can be shown that there exists a linear operator $\alpha:E\to \mathcal{G}^\infty(U)$ such that $\beta\circ\alpha=id_E$  (cf. \cite[Proposition 2.3]{Mujica91}). This clearly parallels the results of \cite{GK} mentioned in Section \ref{sect: Lipschitz free}. Consequently, $E$ is isomorphic to a complemented subspace of $\mathcal{G}^\infty(U)$ whenever $U$ is an open bounded set in $E$. More specifically, we can take $\alpha=D\delta_U(x_0)$ for any $x_0\in U$, where $D\delta_U(x_0)$ denotes the Fr\'{e}chet derivative of $\delta_U$ at $x_0$, which is an operator from $E$ into $\mathcal{G}^\infty(U)$ (recall equation (\ref{eq:frechet-derivative})). When $U=B_E$, we can achieve an isometric version of this \textit{lifting property} by choosing $\alpha=D\delta_U(0)$. \\

In \cite[Corollary 4.12]{Mujica91}, J. Mujica pointed out an alternative description of holomorphic free Banach spaces which is analogous to the identification of Lipschitz free spaces as  \textit{transportation cost spaces} mentioned in the previous section. Given an open subset $U$ of a Banach space $E$, its corresponding holomorphic free Banach space $\mathcal{G}^\infty(U)$ consists of all linear functionals $u\in\mathcal{H}^\infty(U)$ of the form
\begin{equation}\label{eq:functional-holomorphic-free}
 u=\sum_{k=1}^\infty \alpha_k\delta_{x_k},   
\end{equation}
where $(x_k)$ is a sequence of elements of $U$ and $(\alpha_k)$ is an absolutely summable sequence of complex numbers. Furthermore, the norm of an element $u\in \mathcal{G}^\infty(U)$ is given by
$$
\|u\|=\inf \sum_{k=1}^\infty |\alpha_k|,
$$
where the above infimum is taken over all possible representations of $u$ in the form (\ref{eq:functional-holomorphic-free}).
\\

In \cite[Theorem 6.1]{Mujica92}, the author inferred from the simple observation that $\mathcal{H}^\infty(U\times V)$ is isometrically isomorphic to $\mathcal{H}^\infty(U,\,\mathcal{H}^\infty(V))$ and the universal property of projective tensor products (cf. \cite[Proposition 16.13]{BanachSpaceTheory}) that  $\mathcal{G}^\infty(U\times V)=\mathcal{G}^\infty(U)\otimes_{\pi} \mathcal{G}^\infty(V)$. \\

Recall that a sequence $(x_k)_{k=1}^\infty$ of elements of an open subset $U$ of a Banach space $E$ is said to be an \textit{interpolating sequence} if the mapping $$f\in\mathcal{H}^\infty(U)\mapsto (f(x_k) )_{k=1}^\infty\in\ell_\infty$$ is surjective. In \cite[Theorem 5.2]{Mujica912}, J. Mujica relates this notion to holomorphic free Banach spaces, showing that a sequence $(x_k)\subset U$ is interpolating if and only if there exists an operator $T:\mathcal{G}^\infty(U)\to\ell_1$ such that $T\circ S\xi=\xi$ for every $\xi=(\xi_k)\in\ell_1$, where $S$ is the operator defined by
$$
\begin{array}{rclc}
     S:&\ell_1&\longrightarrow &\mathcal{G}^\infty(U) \\
      &(\xi_k)& \longmapsto &\sum_{k=1}^\infty \xi_k\delta_{x_k}.
\end{array}
$$
In particular, it follows that if $U$ has an interpolating sequence, then $\mathcal{G}^\infty(U)$ has a complemented copy of $\ell_1$. This should be compared with the fact that for infinite metric spaces $\mathcal F(M)$ always contains a complemented copy of $\ell_1$ (see Section \ref{sect: Lipschitz free}).

\subsection{Connection between properties of \texorpdfstring{$U$}{Lg} and \texorpdfstring{$\mathcal{G}^\infty(U)$}{Lg}}

We will devote this part of the section to explaining some already known connections between properties of a Banach space $E$ and those corresponding to $\mathcal{G}^\infty(U)$ for $U$ open in $E$.\\ 

The first observation is that the density character of $\mathcal{G}^\infty(U)$ is at most that of the Banach space $E$ containing $U$. By \textit{density character} of a topological space we mean the smallest cardinal for a dense subset of the space. The above claim can be deduced from the fact that $\delta_U(U)$ generates a dense subspace of $\mathcal{G}^\infty(U)$ (cf. \cite[Remark 2.2]{Mujica91}). Note also that for any complex Banach space $E$, by Liouville's theorem $\mathcal H^\infty(E)=\mathbb C$, so $\mathcal G^\infty(E)=\mathbb C$ independently of the density character of $E$.\\

Approximation properties (see Subsection \ref{ss:propertieslipfree}) constitute one of the most studied topics about these connections. 
Given a balanced, bounded and open subset $U$ in $E$, Mujica showed that $E$ has the AP if and only if $\mathcal{G}^\infty(U)$ has the AP \cite[Theorem 5.4]{Mujica91}. In that article, it was pointed out that the aforementioned fact may be helpful for the question of whether $\mathcal{H}^\infty(\mathbb{D})$ has the AP (this problem is formulated in \cite[Problem 1.e.10]{LT1}). \\

 
 In \cite[Proposition 5.7]{Mujica91}, it is established that a Banach space $E$ has the MAP if and only if $\mathcal{G}^\infty(B_E)$ has the MAP. Despite there is some work addressing the analogous question for BAP \cite{Caliskan2}, it is not known whether $E$ has the BAP if and only if $\mathcal{G}^\infty(B_E)$ has the BAP. We thank Ver\'onica Dimant (Buenos Aires) for clarifying this.\\

In line with the preceding results, E.~\c{C}ali\c{s}kan also showed that when $U$ is a bounded balanced open subset of a Banach space $E$, $\mathcal{G}^\infty(U)$ has the CAP if and only if $E$ has the CAP \cite{Caliskan0}. Where a Banach space $E$ is said to have the \textit{compact approximation property} (CAP for short) if for every $\varepsilon>0$ and every $K$ compact set in $E$ there is a compact operator $T:X\to X$ such that $\|Tx-x\|\leq \varepsilon$ for all $x\in K$. \\

Another result worth mentioning is that for any Banach space $E$ the space $\mathcal{G}^\infty(B_E)$ has the Daugavet property, as M. Jung proved in \cite[Corollary 2.3]{Jung} that $\mathcal{H}^\infty(B_E)$ has the Daugavet property and this is inherited by the predual. In particular, $\mathcal{G}^\infty(B_E)$ fails to have the Radon-Nikod\'ym property. As $\mathcal{H}^\infty(B_E)$ cannot have the Radon-Nikod\'ym property either, its predual $\mathcal{G}^\infty(B_E)$ is not an Asplund space (this was previously mentioned in \cite[Example 3.5]{GM}). Nevertheless, we still know very little about general properties of $\mathcal{G}^\infty(U)$.

\subsection{Connection between properties of a holomorphic function and its linearization}

Another issue which must be addressed when we are dealing with free objects is how certain properties of bounded holomorphic functions are \textit{transferred} to their corresponding linear extensions (defined on $\mathcal{G}^\infty(U)$) and vice versa. Let us highlight some examples that illustrate this situation:

\begin{itemize}
    \item Let $E$ and $F$ be Banach spaces and let $U$ be an open subset of $E$. A bounded holomorphic mapping $f:U\to F$ has finite rank if and only if its corresponding operator $\hat{f}:\mathcal{G}^\infty(U)\to F$ has finite rank \cite[Proposition 3.1]{Mujica91}. By a \textit{finite rank} holomorphic mapping $f$ we mean a function such that $\text{span}(f(U))\subset F$ is finite dimensional.

    \item A bounded holomorphic function $f:U\to F$ has a relatively compact range (resp. relatively weakly compact range) if and only if the corresponding operator $\hat{f}:\mathcal{G}^\infty(U)\to F$ is compact (resp. weakly compact) \cite[Proposition 3.4]{Mujica91}. More results in this direction can be found in \cite{J-VR-CS}.

    \item A bounded holomorphic function $f:U\to F$ is \textit{Cohen strongly $p$-summing holomorphic} if and only if its corresponding linearization $\hat{f}:\mathcal{G}^\infty(U)\to F$ is \textit{strongly $p$-summing} (see \cite{J-VSS}). 
\end{itemize}


Consider now $U\subset E$ and $V\subset F$ open subsets of the Banach spaces $E$ and $F$, respectively. Given a holomorphic function $f:U\to V$, we may consider $\delta_V\circ f:U\to \mathcal{G}^\infty(V)$, which is a holomorphic function (as it is a composition of two holomorphic functions). Then, we can define $\overline{f}:\mathcal{G}^\infty(U)\to \mathcal{G}^\infty(V)$ as the unique linear operator which makes commutative the following diagram:
$$
\xymatrix{\mathcal{G}^\infty(U)\ar@{-->}[r]^{\overline{f}}& \mathcal{G}^\infty(V)\\
     U\ar[r]^{f}\ar^{\delta_U}[u]& V\ar^{\delta_V}[u]}
$$

As we have done in the case of Lipschitz free Banach spaces, it is natural to wonder about the relations between $f$ and $\overline{f}$. In \cite[Theorem 5.2]{Mujica92} J. Mujica showed that the \textit{surjectivity} of a holomorphic mapping $f:U\to V$ implies the surjectivity of the operator $\overline{f}:\mathcal{G}^\infty(U)\to \mathcal{G}^\infty(V)$. In contrast, we cannot make an analogous claim about the \textit{injectivity}. For example, let $i:\mathbb{D}\hookrightarrow \mathbb{C}$ be the canonical inclusion of the unit disc $\mathbb{D}$ into the complex plane $\mathbb{C}$. Then $\overline{i}: \mathcal{G}^\infty(\mathbb{D})\to  \mathcal{G}^\infty(\mathbb{C})$ is not injective, since $\mathcal{G}^\infty(\mathbb{C})$ is one-dimensional (as the only elements of $\mathcal{H}^\infty(\mathbb{C})$ are constant functions by Liouville's theorem). \\

Moreover, it can be readily checked that
$f$ is \textit{biholomorphic} if and only if $\overline{f}$ is a \textit{surjective linear isometry}. 
Hence, a natural question motivated by the Lipschitz free spaces case is whether $\mathcal{G}^\infty(U)$ and $\mathcal{G}^\infty(V)$ being isometrically isomorphic implies that their corresponding open sets $U$ and $V$ are biholomorphic. As we shall see below, this is not true in general. \\

\begin{example} Let us denote by $\mathbb{D}^2=\{(z_1,z_2)\in\mathbb{C}^2\::\: |z_1|,|z_2|<1\}$ the complex bidisc. We define $U=3/2\,\mathbb{D}^2$ and $V=U \backslash \{(0,0)\}$. It is clear that the open sets $U$ and $V$ are not biholomorphic (they are not even homeomorphic as topological spaces). 


By Hartogs's extension theorem, given $f\in\mathcal{H}^\infty(V)$, there exists $F\in\mathcal{H}^\infty(U)$ (which is uniquely determined by the Identity Principle) such that $\left.F\right|_V=f$. Thus, we can define an operator $T:\mathcal{H}^\infty(V)\to\mathcal{H}^\infty(U)$ by $Tf=F$. Its linearity can be deduced from the Identity Principle and it is also clear that $T$ is norm preserving and surjective. 

Finally, we are going to show that $T$ is weak$^*$-to-weak$^*$ continuous. By Banach-Dieudonn\'e theorem it is enough to prove that the restriction $\left.T\right|_{B_{\mathcal{H}^\infty(V)}}$ is weak$^*$-to-weak$^*$ continuous. Let $\{f_\alpha\}$ be a weak$^*$-convergent net to some $f$ in $B_{\mathcal{H}^\infty(V)}$. We have to show that $\{T f_\alpha\}$ converges to $Tf$ in $(\mathcal{G}^\infty(U),w^*)$. Since $\{T f_\alpha\}$ is bounded and $\text{span}\{\delta_U(x)\::\: x\in U\}$ is dense in $\mathcal{G}^\infty(U)$, we only need to prove that
$$
F_\alpha(x)=\delta_U(x)(Tf_\alpha)\to \delta_U(x)(Tf)=F(x), \quad \text{ for every } x\in U.
$$
Since $T$ extends functions from $V$ to $U$, the only remaining obstacle is satisfying
$
F_\alpha(0)\to F(0).
$
This can be done by using Cauchy's integral formula. As $F$ is holomorphic on $U$ and $\overline{\mathbb{D}}\subset U$, we have that

$$
F(0,0)=\frac{1}{(2\pi i)^2}\int_{\partial\mathbb{D}}\int_{\partial\mathbb{D}} \frac{F(w_1,w_2)}{w_1w_2}dw_1dw_2=\frac{1}{(2\pi i)^2}\int_{\partial\mathbb{D}}\int_{\partial\mathbb{D}} \frac{f(w_1,w_2)}{w_1w_2}dw_1dw_2,
$$
and hence
$$
    |F_\alpha(0,0)-F(0,0)|  \leq  \|f_\alpha-f\|_{\partial\mathbb{D}\times\partial\mathbb{D},\infty} \to 0.
$$
\end{example}

\begin{rem}
It seems to be an open question whether Banach spaces $E$ and $F$ with holomorphic free Banach spaces $\mathcal{G}^\infty(B_E)$, $\mathcal{G}^\infty(B_F)$ being isometrically isomorphic implies that the unit balls $B_E, \: B_F$ are biholomorphically equivalent. Since W. Kaup and H. Upmeier \cite{KU} (see also \cite{Arazy}) showed that two complex Banach spaces are isometrically isomorphic if and only if their open unit balls are biholomorphically equivalent, the preceding question could be rephrased as follows: are there non-isometrically isomorphic Banach spaces $E,F$ such that $\mathcal{G}^\infty(B_E)$, $\mathcal{G}^\infty(B_F)$ are isometric?
\end{rem}

Finally, we would like to point out some other linearization theorems for certain classes of holomorphic functions which have been considered in the literature. For instance, in \cite{Mazet} (see also \cite{Mujica-Nachbin}) it is shown that for every open subset $U$ of a locally convex space $E$, there exist a complete locally convex space $\mathcal{G}(U)$ and a holomorphic mapping $\delta_U:U\to \mathcal{G}(U)$ such that for each complete locally convex space $F$ and each holomorphic function $f:U\to F$, there is a unique linear mapping $\hat{f}:\mathcal{G}(U)\to F$ such that $\hat{f}\circ \delta_U=f$. Approximation properties of these spaces have been investigated in \cite{Boyd,Caliskan}. There exists an analogous construction for the more restrictive class of \textit{holomorphic functions of bounded type} which can be found in \cite{GGM} and \cite{Mujica92}. A similar study has also been carried out for \textit{weighted holomorphic functions} in \cite{GuBa, GuBa2}. More recently, a \textit{free holomorphic Lipschitz space} has been constructed in \cite{ADGM}. This allows the extension of holomorphic Lipschitz maps to linear operators combining somehow the approaches in this and the previous sections.

\section{Free Banach lattice}\label{sec:fbl}

Let us consider now the category of Banach lattices with lattice homomorphisms and the larger category of Banach spaces with bounded linear operators. Given a Banach space $E$, we can define the \textit{free Banach lattice generated by $E$}. This will be a Banach lattice $\fbl[E]$ together with a linear isometric embedding $\delta_E:E\rightarrow \fbl[E]$, such that for every Banach lattice $X$ and every linear and bounded operator $T:E\rightarrow X$ there is a unique lattice homomorphism $\hat T$ making the following diagram commutative:
$$
\xymatrix{\fbl[E]\ar@{-->}[rd]^{\hat{T}}&\\
     E\ar[r]^{T}\ar^{\delta_E}[u]& X}
$$
Moreover, $\|\hat{T}\|=\|T\|$.\\

\subsection{Construction} 
There exists an explicit construction of $\fbl[E]$ provided by A. Avil\'es, J. Rodr\'iguez and P. Tradacete in \cite{ART}. Namely, $\fbl[E]$ can be represented as a sublattice of the space of positively homogeneous functions from the dual space $E^*$ to $\R$, denoted by $H[E]$, which is a vector lattice when endowed with the pointwise order and lattice operations. It suffices to consider the functional
\begin{equation}\label{equa: norm H1E}
	\nFBLE{f}:=\sup \cbr{\sum_{i=1}^m \abs[0]{f(x_i^*)}:m\in \N, (x_i^*)_{i=1}^m\subset E^*, \sup_{x\in B_E} \sum_{i=1}^m \abs[0]{x_i^*(x)}\leq 1}
\end{equation}
and define the set $H_1[E]:=\cbr{f\in H[E]: \nFBLE{f}<\infty}$, which is a Banach lattice with the norm $\nFBLE{\cdot}$. Note that the evaluation functions $\delta_x:E^*\rightarrow \R$ with $x\in E$, given by $\delta_x(x^*)=x^*(x)$ for $x^*\in E^*$, define a linear isometry
\begin{equation*}
    \fullfunction{\delta_E}{E}{H_1[E]}{x}{\delta_E(x)=\delta_x}
\end{equation*}
It can be checked that the (not necessarily closed) sublattice generated by $\delta_E(E)$ in $H[E]$, denoted by $\FVL[E]$, is the \textit{free vector lattice over $E$ (as a vector space)}, or equivalently, it satisfies the following universal property: given any vector lattice $X$ and any linear operator $T:E\rightarrow X$, there exists a unique lattice homomorphism $\hat{T}:\FVL[E]\rightarrow X$ such that $\hat{T}\circ \delta_E=T$. Note that free vector lattices were already considered in \cite{Baker, Bleier} and their construction follows a standard process. That this vector lattice actually support a lattice norm is less evident. \\

More precisely, the closure of $\FVL[E]$ under the norm $\nFBLE{\cdot}$ will provide the desired representation of $\fbl[E]$. Indeed, the process of extension of linear and bounded operators $T:E\rightarrow X$, $X$ being a Banach lattice, to lattice homomorphism $\hat{T}:\fbl[E]\rightarrow X$ works as follows. First, consider functions $f\in \FVL[E]$, which can be obtained from some evaluation functions $(\delta_{x_i})_{i=1}^n \subset \delta_E(E)$ by means of a finite number of linear and lattice operations, i.e., they can be represented as \textit{lattice-linear expressions} $F(\delta_{x_1},\ldots,\delta_{x_n})$ depending on the $(\delta_{x_i})_{i=1}^n$. Then, defining $\hat{T}f=F(Tx_1,\ldots,Tx_n)$ and using lattice-linear function calculus (cf. \cite[Section 1.d]{LT2}) it can be checked that this expression does not depend on the choice of $F$, so the operator $\hat{T}:\FVL[E]\rightarrow X$ is a well-defined lattice homomorphism. Now, the last step in the proof of \cite[Theorem 2.5]{ART} is to show that $\hat{T}$ is bounded when $\FVL[E]$ is endowed with the norm $\nFBLE{\cdot}$, so that it can be uniquely extended to its closure, $\fbl[E]$.\\

The free Banach lattice over a Banach space $E$ can be visualized as a particular case of a more abstract family of free objects generated by Banach spaces in some restrained classes of Banach lattices, the \textit{free $\Dcal$-convex Banach lattices over a Banach space $E$}, introduced in~\cite{JLTTT}. A \textit{convexity condition} is a triple $\Dcal=(\Gcal, M, \vartheta)$ consisting of the following elements. $\Gcal$ is a nonempty subset of $\bigcup_{m=1}^{\infty} \Hcal_m^{>0}$, where $\Hcal_m$ denotes the set of all continuous positively homogeneous functions from $\R^m$ to $\R$ and $\Hcal_m^{>0}$ is the subset of the monotone functions of $\Hcal_m$, i.e., the functions $h\in \Hcal_m$ such that $h(t_1,\ldots,t_m)\leq h(s_1,\ldots,s_m)$ whenever $0\leq t_i\leq s_i$ for every $i=1,\ldots,m$. On the other hand, $M:\Gcal\rightarrow [1,\infty)$ and $\vartheta:\Gcal \rightarrow \cbr{0,1}$ are just arbitrary functions. Given a convexity condition $\Dcal=(\Gcal, M, \vartheta)$, we say that a seminorm $\nu$ over a vector lattice $X$ is \textit{$\Dcal$-convex} if for every $m\in \N$ and $g\in \Gcal\cap \Hcal_m^{>0}$, the inequality
\begin{equation*}
    \nu(g(x_1,\ldots,x_m))\leq M(g) g(\nu(x_1),\ldots,\nu(x_m))
\end{equation*}
holds for every $x_1,\ldots,x_m\in X_+$ when $\vartheta(g)=1$, and for every pairwise disjoint elements $x_1,\ldots,x_m\in X_+$ when $\vartheta(g)=0$, where the expression $g(x_1,\ldots,x_m)$ must be understood in terms of the lattice-linear function calculus. A Banach lattice $X$ is said to be \textit{$\Dcal$-convex} if its norm is $\Dcal$-convex.\\

Some examples of convexity conditions are $p$-convexity and upper $p$-estimates: Given $1\leq p\leq \infty$, we say that a Banach lattice $X$ is \textit{$p$-convex} (respectively, satisfies an \textit{upper $p$-estimate}) if there exists a constant $M\geq 1$ such that for any finite collection of elements (respectively, pairwise disjoint elements) $x_1,\ldots,x_n\in X$ the inequality
\begin{equation*}
	\norm[3]{\intoo[3]{\sum_{i=1}^n \abs{x_i}^p}^{\frac{1}{p}}}\leq M \intoo[3]{\sum_{i=1}^n \norm{x_i}^p}^{\frac{1}{p}}
\end{equation*}
holds, making the obvious modifications in the case $p=\infty$. We call $p$-convexity (respectively, upper $p$-estimate) constant of $X$ to the least constant $M$ satisfying the above inequality. \\

We can equip the sublattice $\FVL[E]\subset H[E]$ previously defined with the following norm:
\begin{equation*}
    \norm{f}_{\fbl^{\Dcal}[E]}:=\sup\cbr{\nu(f): \nu(\delta_x)\leq \norm{x}\,\,\forall x\in E}
\end{equation*}
for every $f\in \FVL[E]$, where the supremum is taken over all possible $\Dcal$-convex seminorms $\nu$ over $\FVL[E]$. Then, the completion of $\FVL[E]$ under the norm $\norm{\cdot}_{\fbl^{\Dcal}[E]}$, denoted $\fbl^{\Dcal}[E]$, is the free $\Dcal$-convex Banach lattice over $E$ \cite[Theorem 3.3]{JLTTT}, i.e., it is a $\Dcal$-convex Banach lattice that satisfies the following universal property: for every $\Dcal$-convex Banach lattice $X$ and every linear and bounded operator $T:E\rightarrow X$ there exists a unique lattice homomorphism $\hat T$ such that $\hat{T}\circ \delta_E=T$. Moreover, $\norm[0]{\hat{T}}=\norm{T}$. This procedure provides a family of free objects which includes the free Banach lattice over a Banach space $E$ as a particular case. It suffices to observe that, as a consequence of the triangular inequality, being a Banach lattice is equivalent to being 1-convex with constant 1, which is already a convexity condition. Another relevant case is the \textit{free $p$-convex Banach lattice over $E$}, $\fbp[E]$. The existence of this object can be showed in an abstract way by applying the previous procedure with $p$-convexity as our convexity condition, but there also exists an explicit construction that can be obtained by generalizing the norm in \eqref{equa: norm H1E} to define
\begin{equation*}
	\norm{f}_{\fbp[E]}:=\sup \cbr{\intoo[3]{\sum_{i=1}^m \abs[0]{f(x_i^*)}^p}^{\frac{1}{p}}: m\in \N, (x_i^*)_{i=1}^m\subset E^*, \sup_{x\in B_E} \intoo[3]{\sum_{i=1}^m \abs[0]{x_i^*(x)}^p}^{\frac{1}{p}}\leq 1}
\end{equation*}
when $p<\infty$, and
\begin{equation*}
	\norm{f}_{\FBLi[E]}:= \sup_{x^*\in B_{E^*}} \abs[0]{f(x^*)}
\end{equation*}
when $p=\infty$ (see \cite[Section 6]{JLTTT} and \cite[Section 2]{OTTT}). \\

It is worth mentioning that, even though these previous constructions have been developed in the setting of real Banach spaces and Banach lattices, they can be adapted to the complex setting \cite{dHT}. For instance, given a complex Banach space $E$ we can obtain the \textit{free complex Banach lattice over $E$}, $\fbl_{\C}[E]$, through the complexification of an appropriate renorming of a real free Banach lattice. To do so, we consider $E$ as a real Banach space, denoting it by $E_{\R}$, and renorm $\fbl[E_{\R}]$ with
\begin{equation*}
	\norm{f}_{\fbl_{\C}[E]}:=\sup \cbr{\sum_{i=1}^m \abs[0]{f(\mathfrak{Re} z_i^*)}:m\in \N, (z_i^*)_{i=1}^m\subset E^*, \sup_{z\in B_E} \sum_{i=1}^m \abs[0]{z_i^*(z)}\leq 1},
\end{equation*}
which is equivalent to the original norm $\nFBL{\cdot}{E_{\R}}$. Then, the complexified Banach lattice $\fbl_{\C}[E]:=\fbl[E_{\R}]\oplus i\fbl[E_{\R}]$, endowed with the modulus defined pointwise by
\begin{equation*}
    \abs{f}(x^*)=\sqrt{f_1(x^*)^2+f_2(x^*)^2}=\abs{f(x^*)},  x^*\in (E_{\R})^*
\end{equation*}
for every $f=f_1+if_2\in \fbl_{\C}[E]$ and the norm
\begin{equation*}
	\norm{f}_{\fbl_{\C}[E]}:=\sup \cbr{\sum_{i=1}^m \abs[0]{f(\mathfrak{Re} z_i^*)}:m\in \N, (z_i^*)_{i=1}^m\subset E^*, \sup_{z\in B_E} \sum_{i=1}^m \abs[0]{z_i^*(z)}\leq 1},
\end{equation*}
together with the $\C$-linear isometric embedding $\delta_E:E\rightarrow \fbl_{\C}[E]$ given by
$$ \delta_E (z)= \delta_{E_{\R}}(z)-i\delta_{E_{\R}}(iz), $$
plays the role of the free complex Banach lattice over the complex Banach space $E$ \cite[Theorem 3.3]{dHT}.\\

Finally, note that free Banach lattices can be constructed not only over Banach spaces, but over some other categories too. Consider for example the \textit{free Banach lattice generated by a set $A$} proposed by B. de Patger and A. W. Wickstead in \cite{dePW}. It is a Banach lattice $\fbl(A)$ with a bounded map $\iota:A\rightarrow \fbl(A)$ such that for every Banach lattice $X$ and every bounded map $T:A\rightarrow X$ there exists a unique lattice homomorphism $\hat{T}: \fbl(A)\rightarrow X$ with norm $\norm[0]{\hat{T}}=\sup\cbr{\norm{Ta}:a\in A}$ that makes the following diagram commute:
$$
\xymatrix{\fbl(A)\ar@{-->}[rd]^{\hat{T}}&\\
     A\ar[r]^{T}\ar^{\iota}[u]& X}
$$
The free Banach lattice over a set happens to be a particular case of the free Banach lattice over a Banach space, since $\fbl(A)$ is lattice isometric to $\fbl[\ell_1(A)]$ \cite[Corollary 2.9]{ART}. Note that $\ell_1(A)$ can be seen as the free Banach space over the set $A$, because it satisfies the universal property of uniquely extending every bounded map between $A$ and any Banach space $X$ to a bounded linear operator from $\ell_1(A)$ to $X$ by just defining the operator coordinatewise. Thus, the free Banach lattice over a certain set can be interpreted as the free Banach lattice generated by the free Banach space over this set, providing an example of the interaction of free objects in different categories.\\

Free Banach lattices generated by a lattice can also be considered. Recall that a \textit{lattice} $\Lbb$ is a partial ordered set that admits two operations of supremum ($\vee$) and infimum ($\wedge$), and a \textit{lattice homomorphism} in this setting is a map between lattices that preserve these two operations. In \cite{AR1}, A. Avil\'es and J. D. Rodr\'iguez-Abell\'an showed the existence of the \textit{free Banach lattice generated by a lattice $\Lbb$}, which is a Banach lattice $\fbl\langle \Lbb \rangle$ together with a bounded lattice homomorphism $\phi:\Lbb\rightarrow \fbl\langle \Lbb \rangle$ such that for every Banach lattice $X$ and every bounded lattice homomorphism $T:\Lbb \rightarrow X$ there exists a unique Banach lattice homomorphism $\hat{T}: \fbl\langle \Lbb \rangle\rightarrow X$ such that $\hat{T}\circ \phi=T$ and $\norm[0]{\hat{T}}=\sup\cbr{\norm{Tx}: x\in \Lbb}$. They also provided an explicit construction of $\fbl\langle \Lbb \rangle$, inspired by those of $\fbl(A)$ and $\fbl[E]$. Further developments about this notion can be found in \cite{AMRR,AMRR2,AR}\\

\subsection{Properties of $\fbl[E]$}

Let us now focus on some of the properties of $\fbl[E]$. We start by considering the lattice homomorphisms from $\fbl[E]$ to $\R$, which is the same as studying the atoms of $\fbl[E]^*$. These elements are precisely the extensions as lattice homomorphisms of the functionals $x^*\in E^*$, denoted by $\widehat{x^*}$ \cite[Corollary 2.7]{ART}:

\begin{lem}\label{lemm: atoms in FBL*}
    A functional $\phi\in \fbl[E]^*$ is a lattice homomorphism if and only if there exists some $x^*\in E^*$ such that $\phi=\widehat{x^*}$, i.e., $\phi(f)=f(x^*)$ for every $f\in\fbl[E]$.
\end{lem}

Note that, in particular, $\fbl[E]^*$ is a Banach lattice with many atoms, in contrast with other dual Banach lattices such as $L_p[0,1]$, that has no atoms at all.\\

Another interesting fact about $\fbl[E]$ is that its pairwise disjoint subsets can be at most countable, no matter which Banach space $E$ we consider \cite{APR}. This property is known as the \textit{countable chain condition} (\textit{ccc} for short). In fact, $\fbl[E]$ satisfies a stronger property, the \textit{$\sigma$-bounded chain condition} (\textit{$\sigma$-bcc}). Both properties can be stated formally as follows.

\begin{defn}\label{defi: ccc}
    A Banach lattice $X$ 
    \begin{enumerate}
        \item satisfies the \textit{countable chain condition} (\textit{ccc}) if for every uncountable family $\Fcal\subset X_+$ there are distinct $f,g\in \Fcal$ such that $f\wedge g\neq 0$;
        \item satisfies the \textit{$\sigma$-bounded chain condition} (\textit{$\sigma$-bcc}) if $X_+$ admits a countable decomposition $X_+=\bigcup_{n\geq 2} \Fcal_n$ such that, for every $n$, in every subset $\Gcal\subset \Fcal_n$ of size $n$ there exists two distinct elements $f,g\in \Gcal$ such that $f\wedge g\neq 0$.
    \end{enumerate}
\end{defn}

In \cite{APR}, A. Avil\'es, G. Plebanek and J. D. Rodr\'iguez-Abell\'an showed that $\Ccal_{ph}(B_{E^*},\omega^*)$, the space of positively homogeneous and continuous functions over $B_{E^*}$ endowed with the weak* topology, satisfies the $\sigma$-bounded chain condition. Since every function in $\fbl[E]$ is the uniform limit in $B_{E^*}$ of functions in $\FVL[E]$, which are piecewise affine and thus weak* continuous on $B_{E^*}$, it follows that $\fbl[E]$ is a sublattice of $\Ccal_{ph}(B_{E^*},\omega^*)$. In connection with this, it is shown in \cite[Proposition 2.2]{OTTT} that $\Ccal_{ph}(B_{E^*},\omega^*)$ coincides lattice isometrically with $\FBLi[E]$. \\


The case when the underlying Banach space is actually a Banach lattice is of particular interest. First note that in general, for a Banach lattice $X$, $\fbl[X]$ need not coincide with $X$. Take for instance the case of a finite dimensional Banach lattice $X$, where $\fbl[X]$ is lattice isomorphic to $\Ccal(S_{X^*})$ (here, $S_{X^*}$ denotes the unit sphere). This is a common situation in free objects, as the process of embedding $X$ ``freely'' inside $\fbl[X]$ cannot carry the lattice information of $X$. Nevertheless, some of it is preserved in a certain sense. For example, when $X$ is a Banach lattice there exists a surjective lattice homomorphism $\beta_X:\fbl[X]\rightarrow X$, which can be obtained by applying the universal property of $\fbl[X]$ to the identity operator $id_X:X\rightarrow X$. This fact implies that $X$ is complemented (as a Banach space) inside $\fbl[X]$ by means of the projection $\delta_X\circ \beta_X$, whose kernel is an ideal in $\fbl[X]$.\\ 

Inspired by the work on Lipschitz free spaces of G. Godefroy and N. J. Kalton \cite{GK} (see Section \ref{sect: Lipschitz free}), the authors of \cite{AMRT} studied when this situation can be improved, so that $X$ can be embedded inside $\fbl[X]$ through a lattice homomorphism $\alpha_X:X\rightarrow \fbl[X]$ in such a way that $\beta_X\circ \alpha_X=id_X$. A Banach lattice satisfying the above condition is said to have the \textit{lattice-lifting property}. Some examples of Banach lattices with the lattice-lifting property are projective Banach lattices (see \Cref{defi: projective BL} below), Banach lattices with a lattice structure induced by a 1-unconditional basis and $\fbl[E]$ for any Banach space $E$.\\

A particular aspect in the study of the properties of free Banach lattices is to determine which Banach lattice properties of $\fbl[E]$ encode certain Banach space properties of the generating space $E$. In \cite[Section 9]{OTTT} the authors provide a list of characterizations, some of which are summarized in the following theorem.

\begin{thm}\label{theo: equivalences BS BL}
    Let $E$ be a Banach space. Then:
    \begin{enumerate}
        \item $E$ has finite dimension if and only if $\fbl[E]$ has a strong unit;
        \item $E$ is separable if and only if $\fbl[E]$ has a quasi-interior point;
        \item $E$ contains a complemented subspace isomorphic to $\ell_1$ if and only if $\fbl[E]$ contains a sublattice isomorphic to $\ell_1$;
        \item The identity operator $id_{E^*}:E^*\rightarrow E^*$ is $(p',1)$-summing if and only if $\fbl[E]$ satisfies an upper $p$-estimate, where $1\leq p\leq \infty$ and $\frac{1}{p}+\frac{1}{p'}=1$.
    \end{enumerate}
\end{thm}

 This dictionary between Banach space properties of $E$ and Banach lattice properties of $\fbl[E]$ can also be expanded with some entries related to linear and bounded operators between Banach spaces and the lattice homomorphisms between free Banach lattices generated by them. Given Banach spaces $E$ and $F$ and $T:E\rightarrow F$ a linear and bounded operator, we can consider the lattice homomorphism $\overline{T}=\widehat{(\delta_F\circ T)}: \fbl[E]\rightarrow \fbl[F]$, which makes the following diagram commutative:
$$
\xymatrix{\fbl[E]\ar@{-->}[r]^{\overline{T}}& \fbl[F]\\
     E\ar[r]^{T}\ar^{\delta_E}[u]& F \ar^{\delta_F}[u]}
$$
Then, $T$ and $\overline{T}$ are related in the following way \cite[Section 3.1]{OTTT}:

\begin{thm}\label{theo: properties of Tbar}
    Let $E$ and $F$ be Banach spaces and $T:E\rightarrow F$ and $\overline{T}: \fbl[E]\rightarrow \fbl[F]$ as above. Then:
    \begin{enumerate}
        \item $\overline{T}$ is injective if and only if $T$ is injective;
        \item $\overline{T}$ has dense range if and only if $T$ has dense range;
        \item $\overline{T}$ is onto if and only if $T$ is onto;
        \item $\overline{T}$ is an embedding if and only if $T$ is an embedding and there exists some $C>0$ such that for every $n\in\N$, $\eps>0$ and $S:E\rightarrow \ell_1^n$ there exists an extension $\widetilde{S}:E\rightarrow \ell_1^n$ with $\norm[0]{\widetilde{S}}\leq C(1+\eps)\norm{S}$.
    \end{enumerate}
\end{thm}

Note that the previous result implies that if $E$ and $F$ are isomorphic, then $\fbl[E]$ and $\fbl[F]$ are lattice isomorphic. Nevertheless, the reverse question is still open: namely, suppose that $\fbl[E]$ and $\fbl[F]$ are lattice isomorphic, are $E$ and $F$ necessarily linearly isomorphic? Note that in the complex case, this question has a negative answer (see \cite{dHT}).

\subsection{Applications to Banach lattice theory}
Since they were first proposed in 2015 by B. de Patger and A. W. Wickstead \cite{dePW}, free Banach lattices have been proved to be a very useful tool in Banach lattice theory. Take for instance the case of projectivity in Banach lattices, a concept that was also introduced in \cite{dePW}.

\begin{defn}\label{defi: projective BL}
    A Banach lattice $X$ is called \textit{projective} if there exists some $\lambda\geq 1$ such that whenever $Y$ is a Banach lattice, $J\subset Y$ is a closed ideal and $Q: Y\rightarrow Y/J$ is the quotient map, then for every lattice homomorphism $T:X\rightarrow Y/J$ there exists a lattice homomorphism $\Tilde{T}:X\rightarrow Y$ satisfying $Q\circ \Tilde{T}=T$ and $\norm[0]{\Tilde{T}}\leq \lambda \norm{T}$.
\end{defn}

It turns out that free Banach lattices characterize projective Banach lattices in the following sense \cite[Theorem 10.3]{dePW}:

\begin{thm}\label{theo: characterization projectivity}
    A Banach lattice $X$ is projective if and only if for every $\eps>0$ there is a set $A$ and a sublattice $H$ of $\fbl(A)$ such that there exists a lattice isomorphism $\psi: X\rightarrow H$ with $\norm{\psi}$, $\norm[0]{\psi^{-1}}\leq 1+\eps$ and $H$ is lattice complemented in $\fbl(A)$ through a projection of norm at most $1+\eps$.
\end{thm}
Further research concerning the relation between projectivity and free Banach lattices can be found in \cite{AMR2}.\\

Other applications of free Banach lattices to Banach lattice theory can be found in \cite{ART}, where they provide a counterexample to a question raised by J. Diestel on weakly compactly generated Banach lattices. Before formulating the question, let us introduce some key concepts related to the problem. We say that a Banach space $X$ is \textit{weakly compactly generated} (WCG for short) if there exists a weakly compact set $K\subset X$ such that $X=\overline{\spn}(K)$. Analogously, a Banach lattice $X$ is said to be \textit{lattice weakly compactly generated} (LWCG) if there is a weakly compact set $L\subset X$ such that the sublattice generated by $L$ is dense in $X$, i.e., $X=\overline{\latt}(L)$. Clearly, every WCG Banach lattice is LWCG. The question raised by J. Diestel was the converse: is every LWCG Banach lattice WCG? The answer is negative in general, since $\fbl[\ell_2(\Gamma)]$ with $\Gamma$ uncountable is LWCG (take the image of the unit ball of $\ell_2(\Gamma)$ through the canonical embedding $\delta_{\ell_2(\Gamma)}$ as a weakly compact set that generates the whole $\fbl[\ell_2(\Gamma)]$ as a lattice) but not WCG \cite{ART}.\\

An additional open question that has been answered using free Banach lattices comes from \cite{TT}. Recall that a sequence of non-zero vectors $(x_k)$ of a Banach lattice is \textit{bibasic} if there exists a constant $M\geq 1$ such that for every $m\in \N$ and scalars $a_1,\ldots,a_m$ the following inequality, called \textit{bibasic inequality}, is satisfied:
\begin{equation*}
    \norm[3]{\sum_{n=1}^m\abs[3]{\bigvee_{k=1}^n a_k x_k}} \leq M\norm[3]{\sum_{k=1}^m a_kx_k}.
\end{equation*}
It was asked in \cite[Remark 4.4]{TT} whether there exists a subspace of a Banach lattice with no bibasic sequences. In \cite[Theorem 7.5]{OTTT}, T. Oikhberg, M. A. Taylor, P. Tradacete and V. G. Troitsky showed that the canonical copy $\delta_{c_0}(c_0)$ in $\fbl[c_0]$ is a subspace without a bibasic sequence, providing a positive answer to the previous question.\\

Free Banach lattices have also been used to prove the existence of \textit{push-outs} in the category of Banach lattices. In category theory, given three objects $X_0$, $X_1$ and $X_2$ and two morphisms $\alpha_i:X_0\rightarrow X_i$, $i=1,2$, a \textit{push-out} (also known as \textit{amalgamated sum} or \textit{fibered coproduct}) is an object $PO$ together with two morphism $\beta_i:X_i\rightarrow PO$, $i=1,2$, such that the diagram
$$
\xymatrix{X_1\ar@{-->}^{\beta_1}[r]& PO\\
    X_0\ar_{\alpha_2}[r]\ar^{\alpha_1}[u]& X_2 \ar@{-->}_{\beta_2}[u]}
$$
commutes and the following universal property is satisfied: if there exists another object $B$ and two morphisms $\beta'_i:X_i\rightarrow B$, $i=1,2$, such that $\beta'_1\circ \alpha_1=\beta'_2\circ \alpha_2$, then there is a unique $\gamma:PO\rightarrow B$ such that $\beta'_i= \gamma\circ \beta_i$ for $i=1,2$, as the following diagram illustrates.
$$
\xymatrix{ & & B\\
    X_1\ar^{\beta_1}[r] \ar^{\beta'_1}@/^1.5pc/[rru]& PO \ar@{-->}^{\gamma}[ru] &\\
    X_0\ar_{\alpha_2}[r]\ar^{\alpha_1}[u]& X_2 \ar_{\beta_2}[u] \ar_{\beta'_2}@/_1.5pc/[ruu] &}
$$
In the category of Banach lattices with lattice homomorphisms, we can also define an \textit{isometric push-out} as a push-out such that $\norm{\beta_i}\leq 1$, $i=1,2$, and $\norm{\gamma}\leq \max\cbr{\norm{\beta'_1},\norm{\beta'_2}}$ in the universal property. In \cite[Theorem 4.3]{AT} A. Avil\'es a P. Tradacete proved the existence of isometric push-outs using the free Banach lattice generated by a Banach space as an auxiliary tool. Moreover, they proved the following embedding property as a corollary (see also \cite{Tursi}).

\begin{thm}\label{theo: amalgamation embedding of BL}
    Given $T_1:X_0\rightarrow X_1$ and $T_2:X_0\rightarrow X_2$ isometric lattice embeddings between Banach lattices, then there exists a Banach lattice $X_3$ and isometric lattice embeddings $S_i: X_i\rightarrow X_3$, $i=1,2$, such that $S_1\circ T_1=S_2\circ T_2$.
\end{thm}

Free Banach lattices have also been used in \cite{Norm-attaining} to provide the first example of a lattice homomorphism not attaining its norm.

\end{document}